\documentclass[]{interact}

\usepackage{amsmath}
\usepackage{amssymb}
\usepackage{amsfonts}
\usepackage{algorithm}
\usepackage{algorithmic}
\usepackage{mathtools}
\usepackage{xcolor}

\DeclareMathOperator*{\argmax}{argmax}

\newcommand*\dd{\mathop{}\!\mathrm{d}}

\newcommand{\EE}{\mathbb{E}}

\usepackage[natbibapa,nodoi]{apacite}
\renewcommand\bibliographytypesize{\fontsize{10}{12}\selectfont}

\theoremstyle{plain}
\newtheorem{theorem}{Theorem}[section]
\newtheorem{lemma}[theorem]{Lemma}
\newtheorem{corollary}[theorem]{Corollary}

\theoremstyle{definition}
\newtheorem{definition}[theorem]{Definition}

\newtheorem{assumption}{Assumption}
\newtheorem{probstat}{Problem}

\theoremstyle{remark}
\newtheorem{remark}{Remark}

\begin{document}


\title{Analysis of a Generalized Expectation-Maximization Algorithm for Gaussian Mixture Models: A Control Systems Perspective}

\author{
\name{Sarthak Chatterjee\textsuperscript{a}\thanks{CONTACT Sarthak Chatterjee. Email: \texttt{chatts3@rpi.edu}}, Orlando Romero\textsuperscript{b}, and S\'ergio Pequito\textsuperscript{b}}
\affil{\textsuperscript{a}Department of Electrical, Computer, and Systems Engineering, Rensselaer Polytechnic Institute, Troy NY, 12180, USA; \textsuperscript{b}Department of Industrial and Systems Engineering, Rensselaer Polytechnic Institute, Troy NY, 12180, USA}
}

\maketitle

\begin{abstract}
The Expectation-Maximization (EM) algorithm is one of the most popular methods used to solve the problem of parametric distribution-based clustering in unsupervised learning. In this paper, we propose to analyze a generalized EM (GEM) algorithm in the context of Gaussian mixture models, where the maximization step in the EM is replaced by an increasing step. We show that this GEM algorithm can be understood as a linear \mbox{time-invariant} (LTI) system with a feedback nonlinearity. Therefore, we explore some of its convergence properties by leveraging tools from robust control theory. Lastly, we explain how the proposed GEM can be designed, and present a pedagogical example to understand the advantages of the proposed approach.
\end{abstract}


\section{Introduction}

A fundamental problem in unsupervised learning is the problem of \emph{clustering}, where the task in question is to group certain objects of interest into subgroups called \emph{clusters}, such that all objects in a particular cluster share features (in some predefined sense) with each other, but not with objects in other clusters~\citep{Tan:2005:IDM:1095618,Bishop2006}.

The \mbox{Expectation-Maximization} (EM) algorithm~\citep{Dempster1977maximum} is one of the most commonly used methods in parametric \mbox{distribution-based} clustering analysis~\citep{nowak}. Notably, Gaussian mixture models (GMMs) (i.e., a finite convex combination of multivariate Gaussian distributions) have found several applications in real-world problems~\citep{Tan:2005:IDM:1095618,Bishop2006}. In this setting, clustering consists of estimating the parameters in a GMM that maximize its likelihood function (iteratively maximized through the EM algorithm), followed by assigning to each data point the `cluster' corresponding to its most likely multivariate Gaussian distribution in the GMM.



The convergence of the EM algorithm is well-studied in the literature~\citep{wu1983convergence}, particularly in the context of determining the parameters of GMMs~\citep{JordanGMM}. Nonetheless, it is worth analyzing the EM algorithm as a dynamical system, and possibly gain insights that enable us to design more efficient variations of the EM algorithm. Therefore, in~\cite{romero2018convergence}, the authors proposed to change the perspectives on local optimizers and convergence of the EM algorithm by assessing, respectively, the equilibria and asymptotic stability (in the sense of Lyapunov) of a nonlinear dynamical system that represents the standard EM algorithm, through explicit use of discrete-time Lyapunov stability theory.

In this paper, we build upon the recent work in~\cite{romero2018convergence} and propose to analyze a generalized EM (GEM) algorithm~\citep{Dempster1977maximum,neal1998view} in the context of Gaussian mixture models, where the maximization step in the EM is replaced by an increasing step. GEM algorithms have also been used in applications such as computer vision~\citep{fessler1995penalized} and noise estimation in communication channels~\citep{krisjansson2001joint}, and, in general, the study of the EM algorithm and its myriad variants constitute an active area of research~\citep{moon1996expectation,roche2011algorithm}. The main contributions of this work are as follows. First, we show that this GEM algorithm can be understood as a linear \mbox{time-invariant} (LTI) system with a feedback nonlinearity. Secondly, we explore some of its convergence properties by leveraging tools from robust control theory. Lastly, we explain how the proposed GEM can be designed, and present a pedagogical example to understand the advantages of the proposed approach.

\section{Problem Statement}\label{sec:problemstatement}



Let $\theta\in\Theta\subseteq\mathbb{R}^p$ be some vector of unknown (but deterministic) parameters characterizing a distribution of interest, which we seek to infer from a collected dataset $y\in\mathbb{R}^d$ (from now on assumed fixed) and a statistical model composed by a family of joint probability density or mass functions (possibly mixed) $p_\theta(x,y)$ indexed by $\theta\in\Theta$, where 
\begin{equation}
    x\in\mathcal{X} \coloneqq \left\{ x\in\mathbb{R}^n: p_\theta(x,y)>0 \right\}
\end{equation}
is some latent (hidden) random vector. 

The EM algorithm seeks to find a local maximizer of the \emph{incomplete} likelihood function \mbox{$\mathcal{L}:\Theta\to\mathbb{R}$} given by
\begin{equation}
    \mathcal{L}(\theta) \coloneqq p_\theta(y) = \begin{cases}\int_\mathcal{X} p_\theta(x,y)\dd x, &\textnormal{if }x\textnormal{ is continuous,}\\ \sum\limits_{x\in\mathcal{X}} p_\theta(x,y), &\textnormal{if }x\textnormal{ is discrete.}\end{cases}
\end{equation}
The mapping $\theta\mapsto p_\theta(x,y)$ is, naturally, referred to as the \emph{complete} likelihood function. To optimize $\mathcal{L}(\theta)$, the EM algorithm alternates at each iteration $k$ between two steps. First, in the \emph{expectation} step (E-step), we compute $Q(\theta,\theta^{(k)})$, defined through 
\begin{subequations}
\begin{align}
    Q(\theta,\theta') &\coloneqq \EE_{p_{\theta'}(\cdot|y)}[\log p_\theta(\cdot,y)]\label{eq:Qfunction}\\
    &= \int_\mathcal{X} p_{\theta'}(x|y)\log p_\theta(x,y)\dd x,\label{eq:Qintegral}
\end{align}
\end{subequations}
so that $Q(\cdot,\theta^{(k)})$  denotes the expected value of the complete log-likelihood function with respect to $\theta=\theta^{(k)}$. Second, in the \emph{maximization} step (\mbox{M-step}), we maximize $Q(\cdot,\theta^{(k)})$ and update the current estimate as that maximizer.

Before formally stating the EM algorithm, let us make some mild simplifying assumptions that will avoid pathological behavior on the $Q$-function, $Q:\Theta\times\Theta\to\mathbb{R}$.


\begin{assumption}
$\mathcal{X}$ does not depend on $\theta\in\Theta$ and has positive Lebesgue measure.
\label{assump:X}
\end{assumption}

\begin{assumption}
$\mathcal{L}$ is twice continuously differentiable in~$\Theta$.
\label{assump:L}
\end{assumption}
Notice that, from Assumption~\ref{assump:X}, the conditional distribution \mbox{$p_{\theta'}(x|y) = p_{\theta'}(x,y)/p_{\theta'}(y)$} is well defined in $\mathcal{X}$, since $p_\theta(y)>0$ for every $\theta\in\Theta$. Finally, we make the following simplifying assumption, which makes the M-step well defined.

\begin{assumption}
$Q(\cdot,\theta')$ has a unique global maximizer in~$\Theta$.
\label{assum:Q}
\end{assumption}

With all these ingredients and assumptions, we summarize the EM algorithm in Algorithm~\ref{alg:EM}.

\begin{algorithm}
\caption{Expectation-Maximization (EM)}
\textbf{Input:} $y\in\mathbb{R}^d$, $p_\theta(x,y)$, $\theta^{(0)}\in\Theta$.\\
\textbf{Output:} $\hat{\theta}$.
\begin{algorithmic}[1]
\FOR{$k=0,1,2,\ldots$ (until some stopping criterion)} 
    \STATE{\textbf{E-step:} \,\,Compute $Q(\theta,\theta^{(k)})$ }
    \STATE{\textbf{M-step:} Determine $\displaystyle \theta^{(k+1)} = \argmax_{\theta\in\Theta} Q(\theta,\theta^{(k)})$}
\ENDFOR
\RETURN{$\hat{\theta} = $ last computed iteration in $\{\theta^{(k)}\}$.}
\end{algorithmic}
\label{alg:EM}
\end{algorithm}

However, it is to be kept in mind that when we implement the EM algorithm, for most parametric distributions, we do not obtain a closed-form expression for the \mbox{M-step}. As a consequence, to determine a solution (i.e., an approximation) in the \mbox{M-step}, we need to rely on numerical optimization schemes. For instance, we can consider \mbox{first-order} optimization algorithms (e.g., gradient ascent), i.e.,
\begin{equation}
    \theta^{(k+1)} = \theta^{(k)} + \eta\frac{\partial Q(\theta,\theta^{(k)})}{\partial \theta}\Big|_{\theta=\theta^{(k)}},
    \label{eq:Q-GEM}
\end{equation}
for some $\eta>0$. Notice that this could constitute a problem by itself since \mbox{first-order} algorithms are known to have slow convergence rates that get aggravated by the increase in the dimension of the search space. Furthermore, any variant of Algorithm~\ref{alg:EM} that does not explicitly maximize $Q(\cdot,\theta^{(k)})$ at the M-step, but instead is such that $Q(\theta^{(k+1)},\theta^{(k)}) > Q(\theta^{(k)},\theta^{(k)})$ is  referred to as a \emph{generalized} EM (GEM) algorithm.

As previously mentioned, a particularly important class of models are the Gaussian mixture models (GMMs). In these models, each component of the mixture is given by
\begin{equation}
    p_{\theta_i}(y) = \frac{\alpha_i}{\sqrt{\det(2\pi\Sigma_i)}}e^{-\frac{1}{2}(y-\mu_i)^\mathsf{T}\Sigma_i^{-1}(y-\mu_i)},
\end{equation}
with $i = 1,2,\ldots,K$, $y,\mu_i,\in\mathbb{R}^d$, $\Sigma_i\in\mathbb{R}^{d \times d}$ is positive definite, and $\alpha_i\in [0,1]$ such that \mbox{$\sum_{i=1}^K \alpha_i = 1$}. The vector of unknown parameters $\theta$ lumps together the scalar parameters within $\alpha_i,\mu_i,\Sigma_i$ for $i\in\{1,\ldots,K\}$, as follows:
\begin{equation}
    \theta=\left[\alpha^{\mathsf{T}}, \mu^{\mathsf{T}}, \operatorname{vec}\left[\Sigma\right]^\mathsf{T}\right]^\mathsf{T},
\end{equation}
where $\alpha^{\mathsf{T}}=\left[\alpha_1,\ldots, \alpha_K\right]^{\mathsf{T}}$, $\mu^{\mathsf{T}}=\left[\mu_1, \ldots, \mu_K\right]^{\mathsf{T}}$, and \mbox{$\operatorname{vec}\left[\Sigma\right]=\left[
\operatorname{vec}\left[\Sigma_{1}\right]^\mathsf{T}, \ldots, \operatorname{vec}\left[\Sigma_{K}\right]^\mathsf{T}\right]^\mathsf{T}$, with $\operatorname{vec}(M)$} denoting the vector obtained by stacking the column vectors of~$M$.

In this setting, an alternative is to replace the M-step by~\eqref{eq:Q-GEM}, and we obtain a GEM that is able to recover similar (asymptotic) convergence rates available in the literature~\citep{balakrishnan2017statistical}. Nonetheless, (asymptotic)  convergence rates can be misleading as they do not reflect the practical number of iterations required to converge. Furthermore, as it is clear in the GMM, there are some additional constraints that are implicit and are not necessarily satisfied by~\eqref{eq:Q-GEM} (i.e., $\alpha_1+\ldots + \alpha_K = 1$ and $\Sigma_i \succ 0$ for $i = 1,\ldots,K$).

That said, we need to further understand the transient and the local behavior of the GEM algorithm, for which dynamical systems theory provides us with the proper framework. Subsequently, in this paper, we propose to step away from the dynamics without an explicit control (e.g., the M-step in Algorithm~\ref{alg:EM}), towards one where we can consider an additive control, and therefore, study its properties.

In summary, we seek to address the following questions.
\begin{probstat}
\begin{enumerate}
    \item Is it possible to replace the M-step in Algorithm~\ref{alg:EM} by a parameter update step given by
\begin{equation}
    \theta^{(k+1)} = \theta^{(k)} + u^{(k)},
    \label{eq:control}
\end{equation}
where we can design a feedback control law \mbox{$u^{(k)} = \phi(\theta^{(k)})$} to obtain a GEM algorithm?
\item What insights (particularly with respect to design) can such control laws provide us with?
\end{enumerate}
\end{probstat}


\section{Main Results}

In this section, we provide the main result of the paper. Specifically, we show how we can leverage tools from control systems theory to analyze a GEM algorithm as an LTI system connected in feedback with a nonlinearity. Furthermore, we also show how to derive the convergence rate for such an algorithm using tools from robust control. Lastly, we briefly describe how we can look into certain aspects of designing new GEM-like algorithms.

\subsection{GEM Algorithms as LTI Systems with a Feedback Nonlinearity}

We first show how we can leverage tools from dynamical systems and control theory to cast a GEM algorithm into the framework of an LTI system with an interconnected feedback nonlinearity. We begin with the following Lemma that provides us with expressions for the closed-form solution of the problem of estimating the parameters of a GMM using a generalized EM algorithm.

\begin{lemma}[\cite{Dempster1977maximum}]
\label{lemma:lemma1}
Given $K$ possible mixtures in the GMM, and independently and identically distributed (i.i.d.) samples $\{x^{(t)}\}_{t=1}^N$, we can estimate the parameter vector $\theta$ by maximizing the log-likelihood $\mathcal L(\theta)$, that, in the context of a GMM has a closed-form solution given as follows:
\begin{equation}
    \alpha_{j}^{(k+1)}=\frac{1}{N} \sum_{t=1}^{N} h_{j}^{(k)}(t),
\end{equation}
\begin{equation}
    \mu_{j}^{(k+1)}=\frac{1}{\sum_{t=1}^{N} h_{j}^{(k)}(t)} \sum_{t=1}^{N} h_{j}^{(k)}(t) x^{(t)},
\end{equation}
and
\begin{equation}
    \Sigma_{j}^{(k+1)}=\frac{1}{\sum_{t=1}^{N} h_{j}^{(k)}(t)} \sum_{t=1}^{N} h_{j}^{(k)}(t)z_j^{(t),(k+1)}\left(z_j^{(t),(k+1)}\right)^\mathsf{T},
\end{equation}
with
\begin{equation}
    z_j^{(t),(k+1)}=x^{(t)}-\mu_{j}^{(k+1)},
\end{equation}
where the posterior probabilities $h_j^{(k)}$ are given by
\begin{equation}
    h_{j}^{(k)}(t)=\frac{\alpha_{j}^{(k)} p\left(x^{(t)} | \mu_{j}^{(k)}, \Sigma_{j}^{(k)}\right)}{\sum_{i=1}^{K} \alpha_{i}^{(k)} p\left(x^{(t)} | \mu_{i}^{(k)}, \Sigma_{i}^{(k)}\right)}.
\end{equation}
\end{lemma}

With the above closed-form solution, if we, instead, consider a `shifted-update' of the covariance as 
\begin{equation}
\label{eq:shift_cov}
\Sigma_{j}^{(k+1)}=\frac{1}{\sum_{t=1}^{N} h_{j}^{(k)}(t)} \sum_{t=1}^{N} h_{j}^{(k)}(t)z_j^{(t),(k)}\left(z_j^{(t),(k)}\right)^\mathsf{T},
\end{equation}
(i.e., the update of $\Sigma_j^{(k+1)}$ is done with respect to $\mu_j^{(k)}$ instead of $\mu_j^{(k+1)}$), we can summarize in the following Lemma the relationships between the updates of the parameters of the GMM that we aim to estimate, i.e., the mixing weights, the means, and the covariance matrices.
\begin{lemma}
\label{lemma:lemma2}
For the shifted updates of the covariance matrices considered in~\eqref{eq:shift_cov}, the following relations hold:
\begin{equation}\label{eqn:alphaUp}
\alpha^{(k+1)}-\alpha^{(k)}=\left.P_{\alpha^{(k)}} \frac{\partial \mathcal L}{\partial \alpha}\right|_{\alpha=\alpha^{(k)}},
\end{equation}
\begin{equation}\label{eqn:muUp}
\mu_{j}^{(k+1)}-\mu_{j}^{(k)}=\left.P_{\mu_{j}^{(k)}} \frac{\partial \mathcal L}{\partial \mu_{j}}\right|_{\mu_{j}=\mu_{j}^{(k)}},
\end{equation}
and
\begin{equation}\label{eqn:sigmaUp}
\operatorname{vec}\left[\Sigma_{j}^{(k+1)}\right]-\operatorname{vec}\left[\Sigma_{j}^{(k)}\right]=\left.P_{\Sigma_{j}^{(k)}} \frac{\partial \mathcal L}{\partial \operatorname{vec}\left[\Sigma_{j}\right]}\right|_{\Sigma_{j}=\Sigma_{j}^{(k)}},
\end{equation}
with
\begin{align}
    P_{\alpha^{(k)}}&=\frac{1}{N}\left(\operatorname{diag}\left[\alpha_{1}^{(k)}, \ldots, \alpha_{K}^{(k)}\right]-\alpha^{(k)} {\alpha^{(k)}}^\mathsf{T}\right),\\
    P_{\mu_{j}^{(k)}}&=\frac{1}{\sum_{t=1}^{N} h_{j}^{(k)}(t)} \Sigma_{j}^{(k)},\\
    P_{\Sigma_{j}^{(k)}}&=\frac{2}{\sum_{t=1}^{N} h_{j}^{(k)}(t)}\left(\Sigma_{j}^{(k)} \otimes \Sigma_{j}^{(k)}\right),
\end{align}
where $j\in\{1,\ldots, K\}$ denotes the indices of the mixture components, $k$ denotes the iteration number,  and $\otimes$ denotes the Kronecker product. 
\end{lemma}

Therefore, by combining the equations~\eqref{eqn:alphaUp}-\eqref{eqn:sigmaUp} of Lemma~\ref{lemma:lemma2}, we can briefly write the evolution of the parameters as 
\begin{equation}
    \theta^{(k+1)}=\theta^{(k)}+P\left(\theta^{(k)}\right) \nabla \mathcal L(\theta^{(k)}),
\end{equation}
where 
\begin{equation}
    P(\theta)=\operatorname{diag}\left[P_{\alpha}, P_{\mu_{1}}, \ldots, P_{\mu_{K}}, P_{\Sigma_{1}}, \ldots, P_{\Sigma_{K}}\right].
\end{equation}

In this case, the term $\phi(\theta^{(k)})=P\left(\theta^{(k)}\right) \nabla \mathcal L(\theta^{(k)})$ could be understood as a nonlinearity driving the system. Nonetheless, it is not guaranteed that some essential implicit constraints hold, i.e., 
\begin{equation}
\theta^{(k+1)}\in\Theta=\left\{\theta: \sum_{j=1}^{K} \alpha_{j}=1, \Sigma_j=\Sigma_j^{\mathsf{T}}\succ 0 \right\}.
\end{equation}
Therefore, towards incorporating these inter-dependencies, we can consider the `shifted' subspace
\begin{equation}
    \Theta_s=\left\{\theta': \sum_{j=1}^{K} \alpha_{j}=0 \right\},
\end{equation}
where $\theta'=\theta-\theta_0\in \Theta$ for a shift $\theta_0$. Furthermore, let the coordinates of $\theta'$ under the basis $\{e_1,\ldots,e_m\}$ be denoted by $\theta_c$, where the $e_i$-s are canonical (orthonormal) basis vectors and $m$ is the dimension of $\Theta$. Then, $\theta-\theta_0=E\theta_c$, or equivalently, $\theta=E\theta_c+\theta_0$, where $E=[e_1,\ldots,e_m]$. 

Now, notice that $E^\mathsf{T} \theta=E^\mathsf{T} E \theta_{c}+E^\mathsf{T} \theta_{0}$ (or, equivalently, $\theta_{c}=E^\mathsf{T} \theta-E^\mathsf{T} \theta_{0}$), by multiplying both sides by $E^\mathsf{T}$, and noticing that $E^\mathsf{T} E=I$. Thus, $\theta^{\prime}=\theta-\theta_{0}=E^\mathsf{T} E\left(\theta-\theta_{0}\right)=E E^\mathsf{T} \theta^{\prime}\in \Theta$, by observing that $\Theta$ is an open convex set since we only consider local differential properties of the log-likelihood, and consequently, the constraint on positive definiteness of $\Sigma_j$ holds.

Therefore,
\begin{align}
    \theta^{(k+1)}&=E E^\mathsf{T}\left(\theta^{(k)}+P\left(\theta^{(k)}\right) \nabla \mathcal L(\theta^{(k)})\right) \nonumber \\
    &=\theta^{(k)}+ E E^\mathsf{T}P\left(\theta^{(k)}\right) \nabla \mathcal L(\theta^{(k)})
\end{align}
belongs to $\Theta$, which constitutes the parameter update of a GEM algorithm that we shall refer to as \emph{projection-based GEM} (PB-GEM) -- see Algorithm~\ref{alg:GEM}.

\begin{algorithm}
\caption{Projection-Based GEM (PB-GEM)}
\textbf{Input:} $y\in\mathbb{R}^d$, $p_\theta(y)$, $\theta^{(0)}\in\Theta$.\\
\textbf{Output:} $\hat{\theta}$.

\begin{algorithmic}[1]
\FOR{$k=0,1,2,\ldots$ (until some stopping criterion)} 
    \STATE{$\displaystyle \theta^{(k+1)} = \theta^{(k)}+ E E^\mathsf{T}\left.P\left(\theta^{(k)}\right) \frac{\partial \mathcal L}{\partial \theta}\right|_{\theta=\theta^{(k)}}$}
\ENDFOR
\RETURN{$\hat{\theta} = $ last computed iteration in $\{\theta^{(k)}\}$.}
\end{algorithmic}
\label{alg:GEM}
\end{algorithm}

Consequently, the term $\phi(\theta^{(k)})=E E^\mathsf{T}\left.P\left(\theta^{(k)}\right) \frac{\partial \mathcal L}{\partial \theta}\right|_{\theta=\theta^{(k)}}$ can be understood as a nonlinearity driving a linear \mbox{time-invariant} (LTI) system. As such, we can consider \mbox{$\phi(\theta)=\nabla f(\theta)$} to be (locally) Lipschitz and for which there is a (locally) strongly convex function $f$. Before stating the theorem that shows the rate of convergence for the PB-GEM algorithm, we introduce the following preliminary definitions and results.
\begin{definition}[Q-convergence~\citep{jay2001note}]
Given a sequence $\{ \theta^{(k)} \} \to \theta^\star$ with $\theta^{(k)} \neq \theta^\star$ for $k = 0,1,2,\ldots$, the \emph{order} of convergence $\beta$ is a nonnegative number satisfying
\begin{equation}
    \limsup_{k \to \infty} \frac{\| \theta^{(k+1)} - \theta^\star \|}{\| \theta^{(k)} - \theta^\star \|^{\beta}} = \rho < \infty,
\end{equation}
with $\rho$ being the \emph{rate} of convergence.
\end{definition}
\begin{definition}[Sector Integral Quadratic Constraint (IQC) for the gradient map]\label{def:sectorIQC}
For a \emph{strongly convex} function $f$ with strong convexity parameter $\mu_f$, having \emph{Lipschitz continuous} gradients with Lipschitz constant $L$, the gradient map $\nabla f$ satisfies the sector IQC defined by
\begin{equation}\label{eq:sectorIQC}
    \begin{bmatrix}
    \theta - \theta^\star \\ \nabla f(\theta) - \nabla f(\theta^\star)
    \end{bmatrix}^\mathsf{T} \begin{bmatrix} -2 \mu_f L I & (L+\mu_f)I \\ (L+\mu_f)I & -2I \end{bmatrix} \begin{bmatrix}
    \theta - \theta^\star \\ \nabla f(\theta) - \nabla f(\theta^\star)
    \end{bmatrix} \geq 0
\end{equation}
for all $\theta, \theta^\star$.
\end{definition}
\begin{lemma}[A modified version of Theorem 4 in~\cite{lessard2016}]\label{lemma:main_lemma}
Consider a first-order linear optimization scheme represented as the dynamical system
\begin{subequations}\label{eq:dynamic_feedback}
\begin{align}
\xi[k+1] &= A\xi[k] + Bu[k],\label{eq:LTI1}\\
\theta[k] &= C\xi[k] + Du[k],\label{eq:LTI2}\\
u[k] &= \phi(\theta[k]),
\end{align}
\end{subequations}
with nonlinearity $\phi(\theta) = \nabla f(\theta)$. If $\nabla f$ satisfies the sector IQC defined by~\eqref{eq:sectorIQC}, then the linear matrix inequality (LMI)
\begin{equation}
\label{eq:main_LMI}
\begin{bmatrix} A^\mathsf{T}\\B^\mathsf{T} \end{bmatrix} R \begin{bmatrix} A^\mathsf{T}\\B^\mathsf{T} \end{bmatrix}^\mathsf{T} - \begin{bmatrix} \rho^2 R & 0\\ 0 & 0 \end{bmatrix} + \lambda \begin{bmatrix} C & D \\ 0 & I \end{bmatrix}^\mathsf{T} \begin{bmatrix} -2 \mu_f L I & (L+\mu_f)I \\ (L+\mu_f)I & -2I \end{bmatrix} \begin{bmatrix} C & D \\ 0 & I \end{bmatrix} \preceq 0
\end{equation}
is feasible for some $R \succ 0$, $\lambda \geq 0$. Specifically, $\{ \xi[k] \} \to \xi^\star$ with respect to a suitable norm $\| \cdot \|$, with a convergence rate of $\rho$, where $\xi^\star$ is a fixed point of~\eqref{eq:dynamic_feedback} satisfying $\xi^\star = A \xi^\star$.
\end{lemma}

With the above ingredients, we are ready to state our main result concerning the convergence rate of the PB-GEM algorithm, which builds upon tools from robust control theory.

\begin{theorem}\label{thm:first_thm}
Consider a function $f(\theta)$ that is \mbox{$\mu_f$-strongly} convex, has an \mbox{$L$-Lipschitz} gradient, and satisfies \mbox{$\nabla f(\theta) = E E^{\mathsf{T}} P(\theta) \nabla \mathcal{L} (\theta)$}.
Then, $\theta^{(k+1)}=\theta^{(k)}+u^{(k)}$, with \mbox{$u^{(k)} = \nabla f(\theta^{(k)})$}
is a GEM algorithm (i.e., $\{\theta^{(k)} \} \to \theta^\star$, where $\theta^\star$ is the maximum-likelihood estimate) with convergence rate $\rho$  bounded by
\begin{equation}
\label{eq:conv_rate_pb_gem}
\rho \leq \max\{ |1 - \mu_f| , |1 - L| \}.
\end{equation}
\end{theorem}

\begin{proof}
That the Projection-Based GEM algorithm presented in Algorithm~\ref{alg:GEM} indeed constitutes a generalized EM can be shown using an argument similar to one presented in~\cite{SalakhutdinovICML2003}. In particular, if Assumption~\ref{assum:Q} is satisfied for models of the exponential family (a special case being the GMMs considered in this paper), the PB-GEM algorithm evolves in a way such that we have $\mathcal{L}(\theta^{(k+1)}) > \mathcal{L} (\theta^{(k)})$ for all $k \in \mathbb{Z}_+$, provided $\nabla \mathcal{L} (\theta^{(0)}) \neq 0$.
Secondly, we notice that the iterative scheme can be  represented as the LTI system in~\eqref{eq:dynamic_feedback} with a feedback nonlinearity given by \mbox{$\phi(\cdot) = \nabla f(\cdot)$}. 


Due to the regularity of $f$, and the fact that the PB-GEM algorithm can be represented as the dynamical system~\eqref{eq:dynamic_feedback}, with $A = B = C = I$ and $D = 0$, we can invoke the results of Definition~\ref{def:sectorIQC} and Lemma~\ref{lemma:main_lemma} to recover bounds on the convergence rate of the PB-GEM algorithm using the LMI in~\eqref{eq:main_LMI}. Remarkably, due to the general \mbox{block-diagonal} structure of optimization algorithms like gradient ascent, we can then use a `lossless dimensionality reduction argument' and reduce the case of the feasibility of the above LMI to analyze the corresponding semidefinite program for the single-dimensional case without loss of generality~\citep{lessard2016}.

This ascertains the local convergence for the maximum of the function $f$ as long as the following LMI holds 
\begin{equation}
     \begin{bmatrix} (1-\rho^2)R & R \\ R &  R \end{bmatrix} + \lambda \begin{bmatrix} -2 \mu_f L & L+\mu_f \\ L+\mu_f & -2 \end{bmatrix} \preceq 0,
\end{equation}
for some scalar $R > 0$, and $\lambda \geq 0$, where $\rho\in (0,1)$ denotes the convergence rate. Since $R$ is a scalar, we can consider $R = 1$ without loss of generality. This gives us the LMI
\begin{equation}
    \begin{bmatrix} 1-\rho^2 - 2 \mu_f L\lambda & 1 + \lambda(L+\mu_f) \\ 1 + \lambda(L+\mu_f) & 1 - 2\lambda \end{bmatrix} \preceq 0.
\end{equation}
As a consequence, to ensure the negative semidefiniteness of the above matrix, both $1-2\lambda$ (which is present in the bottom right block) and the Schur complement of the bottom right block need to be negative semidefinite. Thus, we have
\begin{subequations}\label{eq:lambda_cond}
\begin{align}
    \lambda & \geq \frac{1}{2},\\
    0 &\geq 1-\rho^2-2 \mu_f L\lambda - \frac{(\lambda(L+\mu_f)+1)^2}{1-2\lambda}.
\end{align}
\end{subequations}
Combining these two, we have
\begin{align}\label{eq:rho_cond}
    \rho^2 &\geq 1 - 2 \mu_f L\lambda - \frac{(\lambda(L+\mu_f)+1)^2}{1-2\lambda},
\end{align}
which yields $\rho \leq \max \{ |1- \mu_f|,|1-  L| \}$.
\end{proof}


Additionally, the transformation matrix $P(\cdot)$ also provides us with valuable insights regarding the rate of convergence of the PB-GEM algorithm. Indeed, differentiating the equation 
\begin{equation}
    F^\mathrm{PB-GEM}(\theta) = \theta + E E^{\mathsf{T}} P(\theta)\nabla\mathcal{L}(\theta),
\end{equation}
we have,
\begin{equation}
     \begin{split}\frac{\partial F^\mathrm{PB-GEM}}{\partial\theta}(\theta) &= I + E E^{\mathsf{T}}\frac{\partial P}{\partial\theta}(\theta)\overline{\overline{\nabla}}\mathcal{L}(\theta) + E E^{\mathsf{T}} P(\theta)S(\theta),\end{split}
\end{equation}
where \mbox{$\frac{\partial P}{\partial\theta} = \begin{bmatrix}\frac{\partial P}{\partial\theta^1} & \ldots   & \frac{\partial P}{\partial\theta^m}\end{bmatrix}$}
with $\theta = \left[ \theta^1,\ldots,\theta^m \right]$,
\begin{equation}
    \overline{\overline{\nabla}}\mathcal{L}(\theta) = \begin{bmatrix} \nabla\mathcal{L}(\theta) & \ldots & 0\\ \vdots & \ddots &  \vdots \\ 0 & \ldots & \nabla\mathcal{L}(\theta)\end{bmatrix},
\end{equation}
and $S(\theta) = \nabla^2 \mathcal{L}(\theta)$ denotes the Hessian matrix of $\mathcal{L}(\cdot)$.

Therefore, near a stationary point $\theta=\theta^\star$ of $\mathcal{L}(\theta)$ over which $P(\theta)$ is bounded, we have
\begin{equation}
    P(\theta) \approx (E E^{\mathsf{T}})^{-1} \left(\frac{\partial F^\mathrm{PB-GEM}}{\partial\theta}(\theta) - I\right)[S(\theta)]^{-1}.
\end{equation}
As a consequence, it follows that, under the aforementioned conditions, the \mbox{projection-based} GEM algorithm exhibits superlinear convergence when $\nabla\mathcal{L}(\theta)$ approaches zero. In particular, the nature of convergence is dictated by the eigenvalues of the matrix $\frac{\partial F^\mathrm{PB-GEM}}{\partial\theta}(\theta)$. If the eigenvalues are near zero, then the transformation matrix scales the EM update step by approximately the scaled negative inverse Hessian, and the EM algorithm behaves like Newton's method. On the other hand, if the eigenvalues are near unity (in absolute value), then the PB-GEM algorithm exhibits \mbox{first-order} convergence.

\subsection{Towards the Design of GEM Algorithms}

We can, therefore, propose to design a GEM algorithm by changing the control law. Nonetheless, we have to be careful with the updates on the different parameters as, implicitly, they possess constraints on the updates. Specifically, we require the $\alpha$-s to sum up to unity, and the $\Sigma$-s to be symmetric positive definite.

Subsequently, in what follows, we focus only on the change of the mean by considering the following weighted function $f_W(\theta)$ such that $f_W(\theta)$ satisfies
\begin{equation}
    \nabla f_W(\theta) = E E^{\mathsf{T}} D P(\theta) \nabla \mathcal{L} (\theta),
\end{equation}
where $D = \text{diag}(I_K,W,I_{d^2 K})$, and with $W\in\mathbb{R}^{dK \times dK}$ being a weight matrix that mixes the different means. In particular, we can consider $W = \text{diag}(\beta_1 I_d,\ldots,\beta_K I_d)$ where $\beta_i>0$ denotes a scaling of the mean similar to a learning rate but applied only on the component rates of the means of the mixture model. Note that we can extend this design step only on the means because the means are the only parameters of the GMMs under consideration that do not have implicit constraints associated with them. This allows us to introduce the following parameter update step
\begin{equation}
    \theta^{(k+1)} = \theta^{(k)}+ E E^\mathsf{T} D P\left(\theta^{(k)}\right) \nabla \mathcal L(\theta^{(k)}),
\end{equation}
for an algorithm which we will refer to as the \emph{weighted} projection-based GEM (\mbox{W-PB-GEM}) algorithm -- see Algorithm~\ref{alg:WGEM}.
\begin{algorithm}
\caption{Weighted Projection-Based GEM (\mbox{W-PB-GEM})}
\textbf{Input:} $y\in\mathbb{R}^d$, $p_\theta(y)$, $\theta^{(0)}\in\Theta$, $W \in \mathbb{R}^{dK \times dK}$.\\
\textbf{Output:} $\hat{\theta}$.
\begin{algorithmic}[1]
\FOR{$k=0,1,2,\ldots$ (until some stopping criterion)} 
    \STATE{$\displaystyle \theta^{(k+1)} = \theta^{(k)}+ E E^\mathsf{T} \left. D P\left(\theta^{(k)}\right) \frac{\partial \mathcal L}{\partial \theta}\right|_{\theta=\theta^{(k)}}$}
\ENDFOR
\RETURN{$\hat{\theta} = $ last computed iteration in $\{\theta^{(k)}\}$.}
\end{algorithmic}
\label{alg:WGEM}
\end{algorithm}

As a result, we have the following corollary on the convergence rate of the \mbox{W-PB-GEM} algorithm.
\begin{corollary}\label{thm:second_thm}
Suppose that there exists a function $f_W(\theta)$ that is~$\mu_f$-strongly convex, has an $L$-Lipschitz gradient, and satisfies $\nabla f_W(\theta) = E E^{\mathsf{T}} D P(\theta) \nabla \mathcal{L} (\theta)$,
where \mbox{$D = \text{diag}(I_K,W,I_{d^2 K})$}, and with \mbox{$W \in \mathbb{R}^{dK \times dK}$} being the matrix of weights that determine the \mbox{component-wise} mixture of the means of the GMM whose parameters are to be estimated.
Then, $\theta^{(k+1)}=\theta^{(k)}+u^{(k)}$, with \mbox{$u^{(k)}=\nabla f_W(\theta^{(k)})$}
is a GEM algorithm (i.e, $\{\theta^{(k)} \} \to \theta^\star$, where $\theta^\star$ is the maximum-likelihood estimate) with convergence rate $\rho$  bounded by
\begin{equation}
\label{eq:conv-w-pb-gem}
\rho \leq \max \{ |1 - \mu_f| , |1 - L| \}.
\end{equation}
\end{corollary}

\begin{remark}
The convergence rates obtained for the \mbox{W-PB-GEM} algorithm are the same as those obtained for the PB-GEM algorithm. It is to be noted, however, that the update equations associated with the $\alpha$-s and the $\Sigma$-s cannot be arbitrarily changed because of the explicit constraints associated with them.
\end{remark}
\begin{remark}
It is worth mentioning here that the convergence rates as obtained in~\eqref{eq:conv_rate_pb_gem} and~\eqref{eq:conv-w-pb-gem} are merely upper bounds, and, unfortunately, do not shed any light on the transient behavior of the \mbox{PB-GEM} or \mbox{W-PB-GEM} algorithm -- see the inset of Figures~\ref{fig:pb-gem-2} and~\ref{fig:w-pb-gem-2}.
\end{remark}

\section{Pedagogical Examples}

In this section, we seek to demonstrate a pedagogical example that shows the efficacy of the methods extended in this paper in identifying the parameters of unknown GMMs. To do this, we first sample $1000$ arbitrary points from a mixture of two Gaussians with the following parameters
\[ \mu_1^\star = \begin{bmatrix} 1\\1 \end{bmatrix} \: \text{and} \: \mu_2^\star = \begin{bmatrix} -1\\-1 \end{bmatrix}, \]
\[ \Sigma_1^\star = \begin{bmatrix} 1 & 0 \\ 0 & 1 \end{bmatrix} \: \text{and} \: \Sigma_2^\star = \begin{bmatrix} 1 & 0 \\ 0 & 1 \end{bmatrix}, \]
and
\[ \alpha^\star = \begin{bmatrix} 0.5 \\ 0.5 \end{bmatrix}. \]
Further, we initialized the algorithms with the following parameters
\[ \mu_1 = \begin{bmatrix} -3 \\ 3 \end{bmatrix}, \mu_2 = \begin{bmatrix} 3 \\ -3 \end{bmatrix}, \]
such that they lie on the line which is orthogonal to the direction defined by $\mu_1^\star$ and $\mu_2^\star$. Additionally, $\Sigma_1$ and $\Sigma_2$ are initialized to be arbitrary positive definite matrices and $\alpha$ is arbitrarily initialized such that \mbox{$\alpha_1 + \alpha_2 = 1$}.

\subsection{The PB-GEM algorithm}

\begin{figure}[ht]
    \centering
    \includegraphics[width=0.8\textwidth]{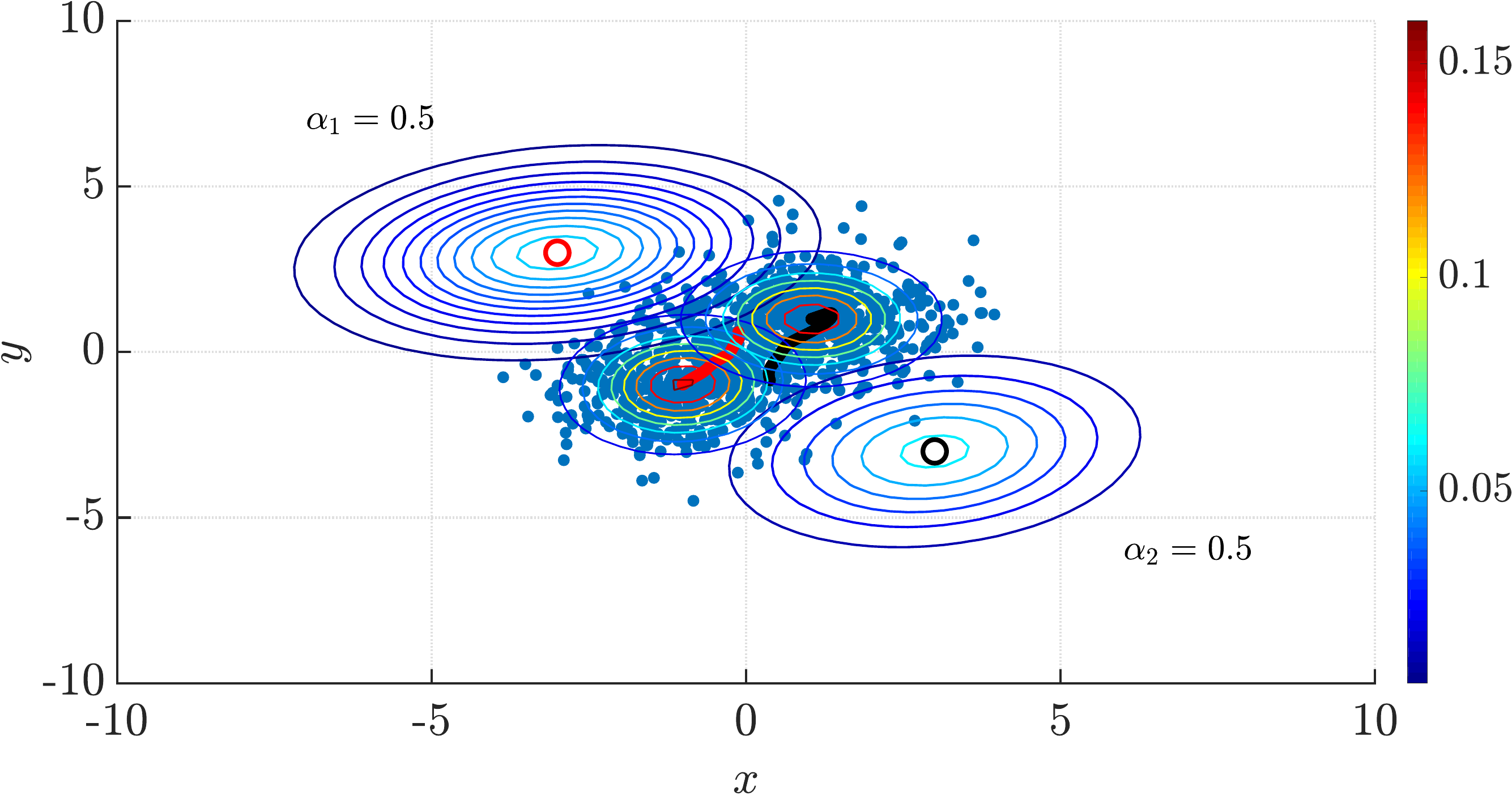}
    \caption{1000 randomly sampled points from a synthetic two-component Gaussian mixture model and the contour plots of the initial and final estimated probability distribution functions using the PB-GEM algorithm. The red and black circles show the initial positions, and the red and black traces show the evolutions of $\mu_1$ and $\mu_2$, respectively.}
    \label{fig:pb-gem-1}
\end{figure}

\begin{figure}[ht]
    \centering
    \includegraphics[width=0.8\textwidth]{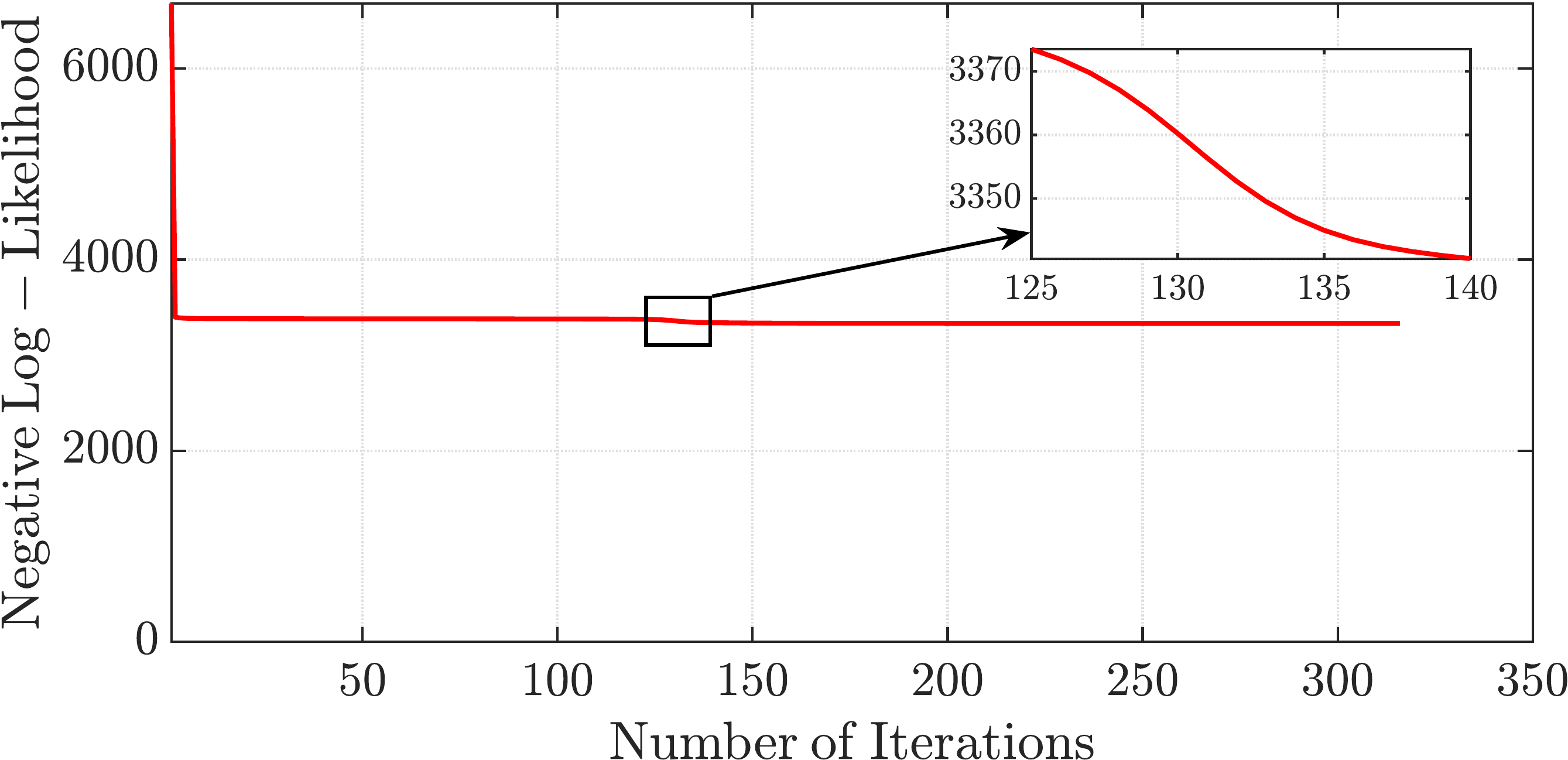}
    \caption{Plot of negative log-likelihood versus the number of iterations for the PB-GEM algorithm. The transient behavior between iterations 125 and 140 is shown in the inset.}
    \label{fig:pb-gem-2}
\end{figure}

We first consider a pedagogical example to demonstrate the performance of the proposed PB-GEM algorithm on estimating the parameters of the synthetic Gaussian mixture
model specified above. The results of running the \mbox{PB-GEM} algorithm to determine the parameters of the above mixture model are shown in Figures~\ref{fig:pb-gem-1} and~\ref{fig:pb-gem-2}. We see that the proposed PB-GEM algorithm is able to successfully determine the parameters of the synthetic GMM from which the points have been sampled. Convergence is obtained in 316 iterations (i.e., to attain a relative change in log-likelihood smaller than $10^{-10}$).

\subsection{The W-PB-GEM algorithm}\label{subsec:subsec}

\begin{figure}[ht]
    \centering
    \includegraphics[width=0.8\textwidth]{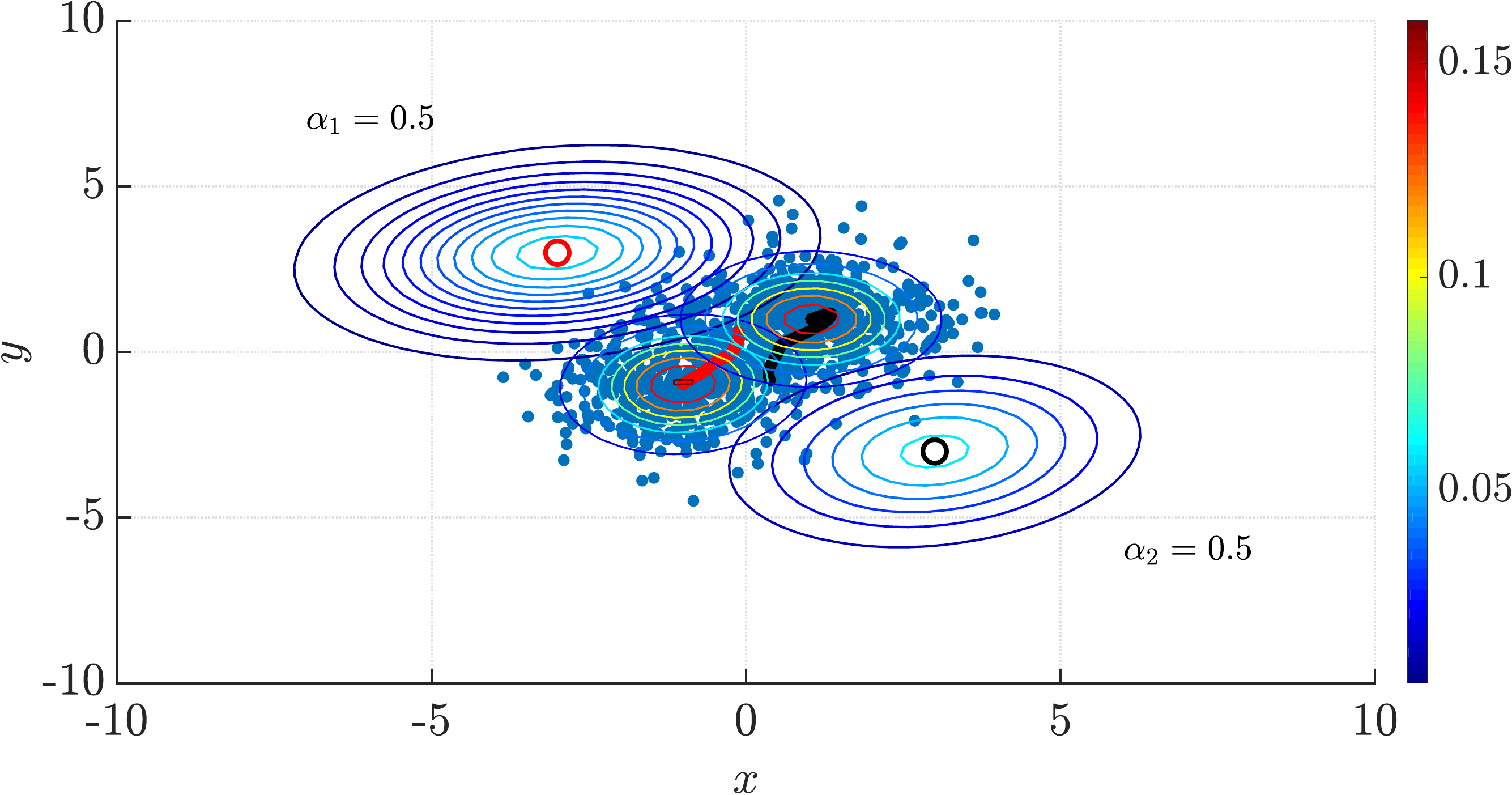}
    \caption{1000 randomly sampled points from a synthetic two-component Gaussian mixture model and the contour plots of the initial and final estimated probability distribution functions using the W-PB-GEM algorithm. The red and black circles show the initial positions, and the red and black traces show the evolutions of $\mu_1$ and $\mu_2$, respectively.}
    \label{fig:w-pb-gem-1}
\end{figure}

Next, we test the performance of the W-PB-GEM algorithm. The matrix of weights $W$ that determine the mixture of proportions during the updates of the means is given by
\[ W = \begin{bmatrix} 0.996 I_2 & 0 \\ 0 & 0.996 I_2 \end{bmatrix}. \]
The results of running the W-PB-GEM algorithm with the same initializations for $\mu_1, \mu_2, \Sigma_1, \Sigma_2$, and $\alpha$ are documented in Figures~\ref{fig:w-pb-gem-1} and~\ref{fig:w-pb-gem-2}. Convergence is obtained in 279 iterations with the same stopping criterion used in the previous section.

\subsection{Multi-Class Classification}
\begin{figure}[ht]
    \centering
    \includegraphics[width=0.8\textwidth]{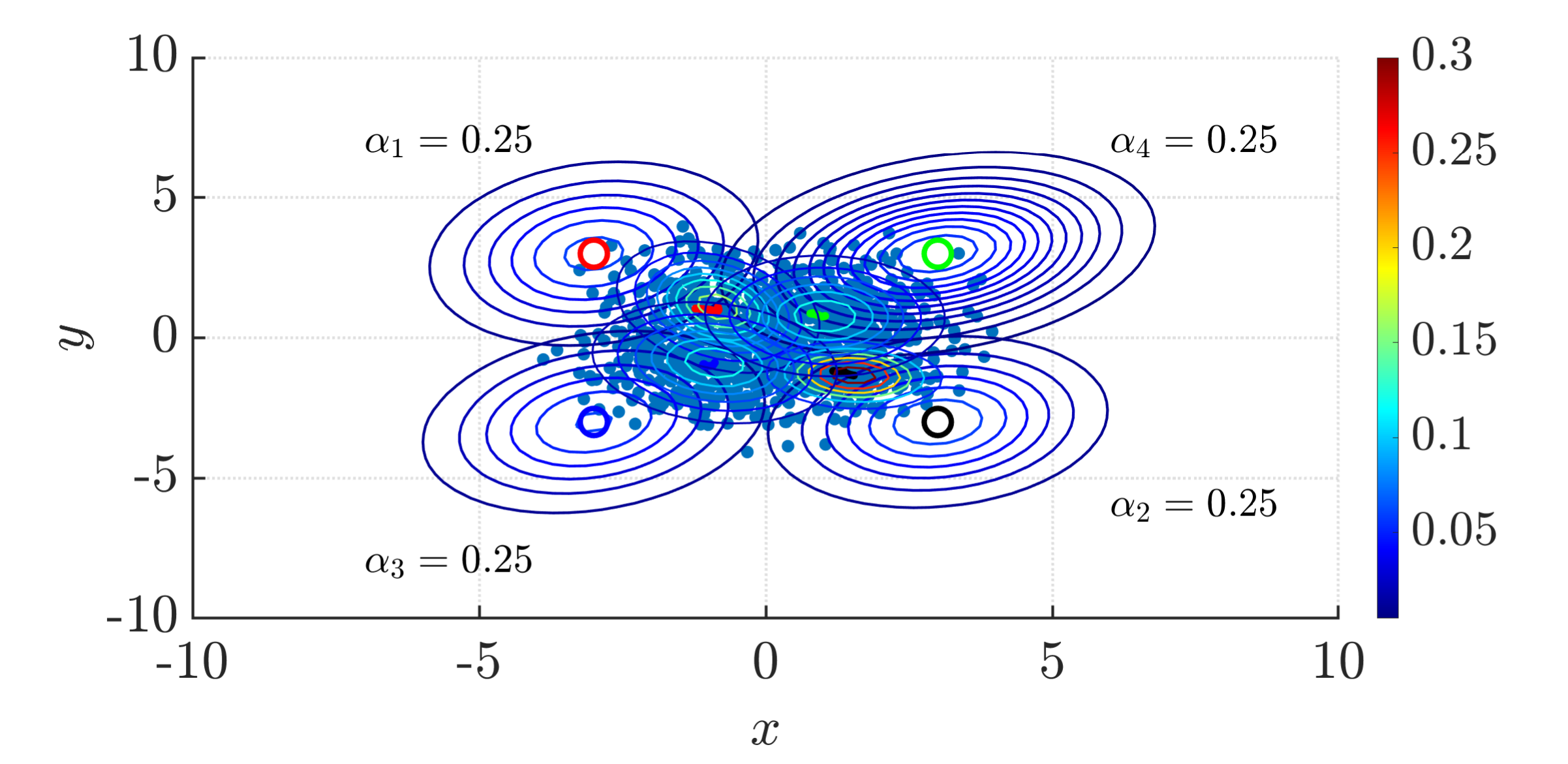}
    \caption{1000 randomly sampled points from a synthetic four-component Gaussian mixture model and the contour plots of the initial and final estimated probability distribution functions using the PB-GEM algorithm. The red, black, blue, and green circles show the initial positions, and the respective traces show the evolutions of $\mu_i, i=1,2,3,4$.}
    \label{fig:multi-pb-gem}
\end{figure}

\begin{figure}[ht]
    \centering
    \includegraphics[width=0.8\textwidth]{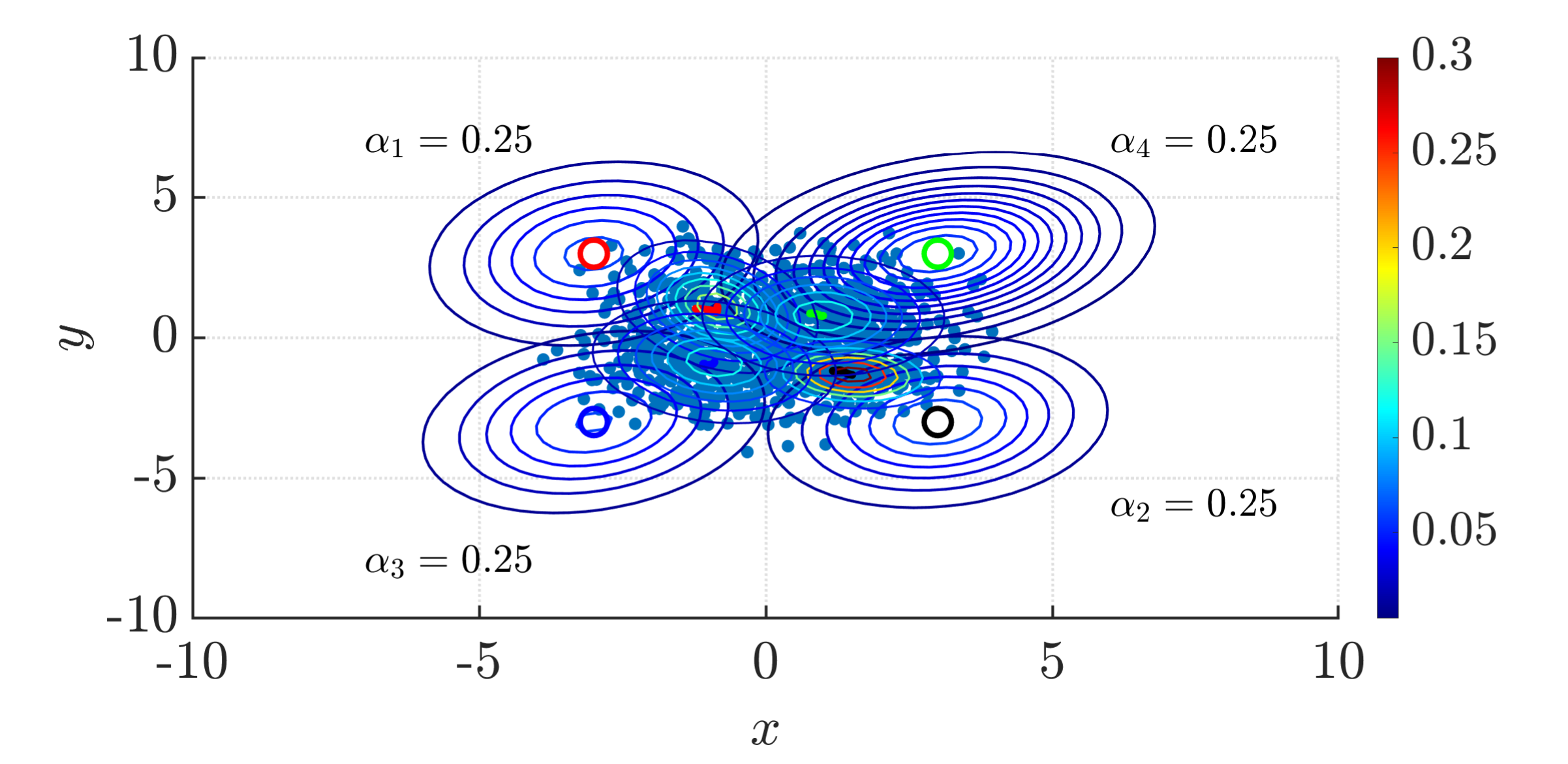}
    \caption{1000 randomly sampled points from a synthetic four-component Gaussian mixture model and the contour plots of the initial and final estimated probability distribution functions using the W-PB-GEM algorithm. The red, black, blue, and green circles show the initial positions, and the respective traces show the evolutions of $\mu_i, i=1,2,3,4$.}
    \label{fig:multi-w-pb-gem}
\end{figure}

In what follows, we also present an illustrative example where both the PB-GEM and the W-PB-GEM algorithms are used in order to identify the parameters of a GMM with more than two Gaussians. We sample $1000$ arbitrary points from a mixture of four Gaussians with the following parameters
\[ \mu_1^\star = \begin{bmatrix} 1\\1 \end{bmatrix}, \mu_2^\star = \begin{bmatrix} -1\\-1 \end{bmatrix}, \mu_3^\star = \begin{bmatrix} 1\\-1 \end{bmatrix}, \: \text{and} \: \mu_4^\star =  \begin{bmatrix} -1\\1 \end{bmatrix}, \]
\[ \Sigma_1^\star = \Sigma_2^\star = \Sigma_3^\star = \Sigma_4^\star = \begin{bmatrix} 1 & 0 \\ 0 & 1 \end{bmatrix}, \]
and
\[ \alpha^\star = \begin{bmatrix} 0.25 \\ 0.25 \\ 0.25 \\ 0.25 \end{bmatrix}. \]
Further, we initialized both the algorithms with the following parameters
\[ \mu_1 = \begin{bmatrix} -3 \\ 3 \end{bmatrix}, \mu_2 = \begin{bmatrix} 3 \\ -3 \end{bmatrix}, \mu_3 = \begin{bmatrix} -3 \\ -3 \end{bmatrix}, \: \text{and} \: \mu_4 = \begin{bmatrix} 3 \\ 3 \end{bmatrix}. \]
Additionally, $\Sigma_1, \Sigma_2, \Sigma_3,$ and $\Sigma_4$ are initialized to be arbitrary positive definite matrices and $\alpha$ is arbitrarily initialized such that $\sum_{i=1}^4 \alpha_i = 1$. The matrix of weights $W$ that determine the mixture of proportions during the updates of the means for the \mbox{W-PB-GEM} algorithm is once again given by
\[ W = \begin{bmatrix} 0.996 I_2 & 0 \\ 0 & 0.996 I_2 \end{bmatrix}. \] The results of running the PB-GEM and W-PB-GEM algorithms for this problem are shown in Figures~\ref{fig:multi-pb-gem} and~\ref{fig:multi-w-pb-gem} respectively. The results are similar to the case on two-class classification. Convergence (i.e., attaining a relative change in log-likelihood smaller than $10^{-10}$) is obtained in $1822$ iterations for the PB-GEM algorithm and in $472$ iterations for the W-PB-GEM algorithm.

\subsection{Discussion of results}

\begin{figure}[ht]
    \centering
    \includegraphics[width=0.8\textwidth]{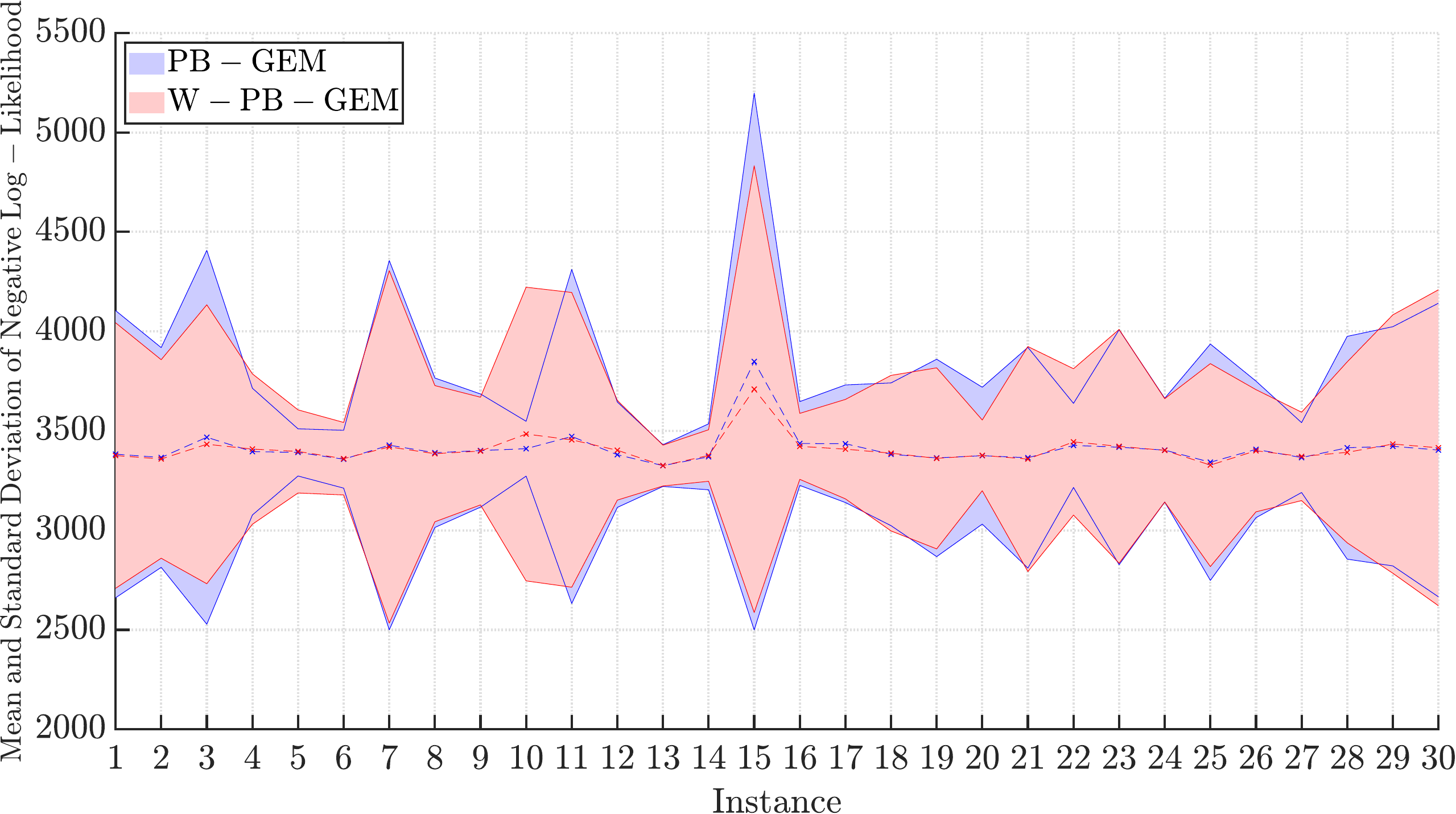}
    \caption{Means and standard deviations of the negative \mbox{log-likelihood} function for 30 instances of the \mbox{PB-GEM} and \mbox{W-PB-GEM} algorithm for the two-class classification example running with the same initial conditions for a particular instance.}
    \label{fig:logLplot}
\end{figure}

The reason why the initial conditions on the means are selected such that they lie on a line orthogonal to the means characterizing the synthetic GMM considered above is to intentionally make the convergence of the \mbox{PB-GEM} and \mbox{W-PB-GEM} algorithms more difficult. We also illustrate in Figure~\ref{fig:logLplot} a comparative study of the \mbox{PB-GEM} and \mbox{W-PB-GEM} algorithms for the two-class example (the matrix $W$ was selected identical to the one in Section~\ref{subsec:subsec}) by plotting the mean and standard deviation of the negative \mbox{log-likelihood} function for 30 instances of both the algorithms when they are initialized with the same set of parameters for a particular instance.
In general, such worst-case initialization conditions are useful in order to gain insights into the transient behaviors of such algorithms.

It is also instructive here to note that for the problem of identifying the parameters of a high-dimensional GMM, the number of iterations to convergence would grow exponentially. In such a case, it would be extremely important to have convergence to the actual parameters in as few iterations as possible, since each iteration would involve a pass over the entire dataset, and when the dataset is large, having a lower number of iterations to convergence would amount to a reduction in the amount of time taken for the estimation of the parameters.

\begin{figure}[ht]
    \centering
    \includegraphics[width=0.8\textwidth]{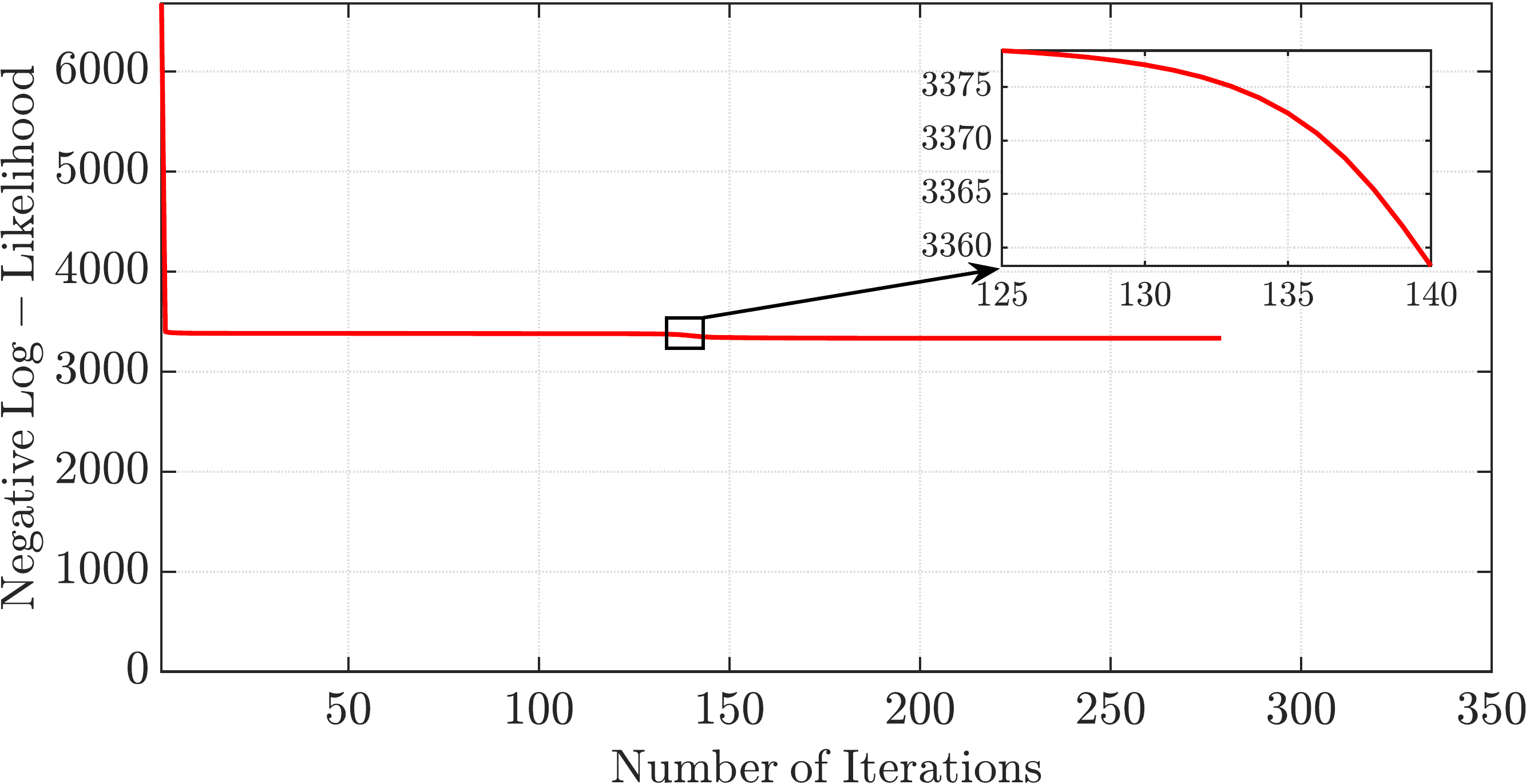}
    \caption{Plot of negative log-likelihood versus the number of iterations for the W-PB-GEM algorithm. The transient behavior between iterations 125 and 140 is shown in the inset.}
    \label{fig:w-pb-gem-2}
\end{figure}

We also reiterate that the convergence rates presented in~\eqref{eq:conv_rate_pb_gem} and~\eqref{eq:conv-w-pb-gem} are merely upper bounds. As demonstrated in the insets of Figures~\ref{fig:pb-gem-2} and~\ref{fig:w-pb-gem-2}, these have no relationship with the transient behaviors of the \mbox{PB-GEM} and \mbox{W-PB-GEM} algorithms. Although in practice we can improve the convergence rates of these algorithms by designing new and more efficient varieties (as detailed in the construction of the \mbox{W-PB-GEM} algorithm), the upper bound of the obtained convergence rates does not change.


\section{Conclusions and Future Work}

In this paper, we analyzed a GEM algorithm to estimate the parameters of GMMs from a dynamical systems perspective. In particular, we showed that this algorithm can be understood as an LTI system connected in feedback with a nonlinearity. The convergence properties of the proposed algorithm are studied by utilizing tools from robust control theory. We also explored the simple design of this class of GEM algorithms and provided evidence using pedagogical examples that it might be possible to improve the transient and the practical convergence of these algorithms despite the fact that they exhibit the same asymptotic convergence rates. Future work will consist of using tools from adaptive systems theory to accelerate practical convergence properties for GEM algorithms. Additionally, fundamental connections exist between the EM algorithm and proximal point methods~\citep{chretien2000kullback,chretien2008algorithms,figueiredo2004lecture} and future work will focus on analyzing proximal interpretations of the EM algorithm using tools from robust control theory~\citep{lessard2016,fazlyab2018analysis}.

\bibliographystyle{apacite}
\bibliography{mybibfile}

\end{document}